
\magnification = \magstep0
\overfullrule = 2pt
\vsize = 525dd
\hsize = 27cc
\topskip = 13dd
\hoffset = 1.8cm
\voffset = 1.8cm


\font\bfbig = cmbx10 scaled \magstep2   
\font\partf = cmbx10 scaled \magstep1	

\font\eightrm = cmr10 scaled 800        
\font\sixrm = cmr7 scaled 850
\font\fiverm = cmr5
\font\eighti = cmmi10 scaled 800
\font\sixi = cmmi7 scaled 850
\font\fivei = cmmi5
\font\eightit = cmti10 scaled 800
\font\eightsy = cmsy10 scaled 800
\font\sixsy = cmsy7 scaled 850
\font\fivesy = cmsy5
\font\eightsl = cmsl10 scaled 800
\font\eighttt = cmtt10 scaled 800
\font\eightbf = cmbx10 scaled 800
\font\sixbf = cmbx7 scaled 850
\font\fivebf = cmbx5

\font\fivefk = eufm5                    
\font\sixfk = eufm7 scaled 850
\font\sevenfk = eufm7
\font\eightfk = eufm10 scaled 800
\font\tenfk = eufm10

\newfam\fkfam
	\textfont\fkfam=\tenfk \scriptfont\fkfam=\sevenfk
		\scriptscriptfont\fkfam=\fivefk

\def\eightpoint{%
	\textfont0=\eightrm \scriptfont0=\sixrm \scriptscriptfont0=\fiverm
		\def\rm{\fam0\eightrm}%
	\textfont1=\eighti  \scriptfont1=\sixi  \scriptscriptfont1=\fivei
		\def\oldstyle{\fam1\eighti}%
	\textfont2=\eightsy \scriptfont2=\sixsy \scriptscriptfont2=\fivesy
	\textfont\itfam=\eightit \def\it{\fam\itfam\eightit}%
	\textfont\slfam=\eightsl \def\sl{\fam\slfam\eightsl}%
	\textfont\ttfam=\eighttt \def\tt{\fam\ttfam\eighttt}%
	\textfont\bffam=\eightbf \scriptfont\bffam=\sixbf
		\scriptscriptfont\bffam=\fivebf \def\bf{\fam\bffam\eightbf}%
	\textfont\fkfam=\eightfk \scriptfont\fkfam=\sixfk
		\scriptscriptfont\fkfam=\fivefk
	\rm}
\skewchar\eighti='177\skewchar\sixi='177\skewchar\eightsy='60\skewchar\sixsy='60

\def\fk{\fam\fkfam}
\def\petit{\eightpoint}%


\def\today{\ifcase\month\or
	January\or February\or March\or April\or May\or June\or
	July\or August\or September\or October\or November\or December\fi
	\space\number\day, \number\year}
\def\newline{\hfil\break}


\newtoks\RunAuthor\RunAuthor={} \newtoks\RunTitle\RunTitle={}
\def\ShortTitle#1#2{\RunAuthor={#1}\RunTitle={#2}}
\headline={\ifnum\pageno=1{\hfil}\else
	\ifodd\pageno {\petit{\the\RunTitle}\hfil\folio}
	\else {\petit\folio\hfil{\the\RunAuthor}} \fi \fi}
\footline={\hfil}
\long\def\Title#1{\hrule\topglue3truecm
	\noindent{\bfbig#1}\vskip12pt\relax}
\long\def\Author#1#2{\noindent{\bf#1}\vskip6pt%
	\noindent{\petit#2}\vskip10pt%
	\noindent{\petit\today}\vskip32pt\relax}
\long\def\Thanks#1#2{{\parindent=20pt\baselineskip=9pt%
	\footnote{\nobreak${}^{#1}$}{\petit#2\par\vskip-9pt}}}
\long\def\Abstract#1{\noindent{\bf Summary.}\enspace#1\bigskip\relax}


\def\PLabel#1{\xdef#1{\nobreak(p.\the\pageno)}}


\newcount\SECNO \SECNO=0
\newcount\SUBSECNO \SUBSECNO=0
\newcount\SUBSUBSECNO \SUBSUBSECNO=0
\def\Part#1{\SECNO=0\SUBSECNO=0\SUBSUBSECNO=0
	\vfill\eject\noindent{\partf Part #1}
	\nobreak\vskip12pt\noindent\kern0pt}
\def\Section#1{\SUBSECNO=0\SUBSUBSECNO=0 \advance\SECNO by 1
	\goodbreak\vskip21pt\leftline{\bf\the\SECNO .\ #1}
	\gdef\Label##1{\xdef##1{\nobreak\the\SECNO}}
	\nobreak\vskip12pt\noindent\kern0pt}
\def\SubSection#1{\SUBSUBSECNO=0 \advance\SUBSECNO by 1
	\goodbreak\vskip21pt\leftline{\it\the\SECNO.\the\SUBSECNO\ #1}
	\gdef\Label##1{\xdef##1{\nobreak\the\SECNO.\the\SUBSECNO}}
	\nobreak\vskip12pt\noindent\kern0pt}
\def\SubSubSection#1{\advance\SUBSUBSECNO by 1
	\goodbreak\vskip21pt\leftline{\rm\the\SECNO.\the\SUBSECNO.\the\SUBSUBSECNO\ #1}
	\gdef\Label##1{\xdef##1{\nobreak\the\SECNO.\the\SUBSECNO.\the\SUBSUBSECNO}}
	\nobreak\vskip12pt\noindent\kern0pt}


\long\def\Definition#1#2{\medbreak\noindent{\bf Definition%
	#1.\enspace}{\it#2}\medbreak\smallskip\relax}
\long\def\Theorem#1#2{\medbreak\noindent{\bf Theorem%
	#1.\enspace}{\it#2}\medbreak\smallskip\relax}
\long\def\Lemma#1#2{\medbreak\noindent{\bf Lemma%
	#1.\enspace}{\it#2}\medbreak\smallskip\relax}
\long\def\Proposition#1{\medbreak\noindent{\bf Proposition.%
	\enspace}{\it#1}\medbreak\smallskip\relax}


\newcount\FOOTNO \FOOTNO=0
\long\def\Footnote#1{\global\advance\FOOTNO by 1
	{\parindent=20pt\baselineskip=9pt%
	\footnote{\nobreak${}^{\the\FOOTNO)}$}{\petit#1\par\vskip-9pt}%
	}\gdef\Label##1{\xdef##1{\nobreak\the\FOOTNO}}}


\newcount\EQNO \EQNO=0
\def\Eqno{\global\advance\EQNO by 1 \eqno(\the\EQNO)%
	\gdef\Label##1{\xdef##1{\nobreak(\the\EQNO)}}}


\newcount\FIGNO \FIGNO=0
\def\Fcaption#1{\global\advance\FIGNO by 1
	{\petit{\bf Fig. \the\FIGNO.~}#1}
	\gdef\Label##1{\xdef##1{\nobreak\the\FIGNO}}}

\newcount\TABNO \TABNO=0
\def\Tcaption#1{\global\advance\TABNO by 1
   {\petit{\bf Table. \the\TABNO.~}#1}
   \gdef\Label##1{\xdef##1{\nobreak\the\TABNO}}}


\newcount\REFNO \REFNO=0
\newbox\REFBOX \setbox\REFBOX=\vbox{}
\def\BegRefs{\setbox\REFBOX\vbox\bgroup
	\parindent18pt\baselineskip9pt\petit}
\def\EndRefs{\par\egroup}
\def\Ref#1{\global\advance\REFNO by 1 \ifnum\REFNO>1\vskip3pt\fi
	\item{\the\REFNO .~}\xdef#1{\nobreak[\the\REFNO]}}
\def\References{\goodbreak\vskip21pt\leftline{\bf References}
	\nobreak\vskip12pt\unvbox\REFBOX\vskip21pt\relax}


\def\N{I\kern-.8ex N}
\def\Z{\raise.72ex\hbox{${}_{\not}$}\kern-.45ex {\rm Z}}
\def\Q{\raise.82ex\hbox{${}_/$}\kern-1.35ex Q} \def\R{I\kern-.8ex R}
\def\C{\raise.87ex\hbox{${}_/$}\kern-1.35ex C} \def\H{I\kern-.8ex H}

\def\D#1#2{{{\partial#1}\over{\partial#2}}}



\Title{On Guichard's nets and Cyclic systems}
\Author{U.~Hertrich-Jeromin\Thanks{\ast}{Partially supported by
	DFG grant He 2490/1-1 and by the MSRI at ETH Z\"urich.},
	E.-H.~Tjaden and M.T.~Z\"urcher}
	{Dept Math \&\ Stat, GANG, University of Massachusetts,
	 Amherst, MA 01003\newline
	Dept Mathematics, SFB288, Technical University Berlin,
	 D-10623 Berlin\newline
	Dept Mathematics, ETH-Zentrum, CH-8092 Z\"urich}
\ShortTitle{U.Jeromin, E.-H.Tjaden, M.T.Z\"urcher}
	{Guichard's nets and Cyclic Systems}

\BegRefs
\Ref\Blaschke W.~Blaschke: {\it Vorlesungen \"uber Differentialgeometrie
	III\/}; Springer, Berlin 1928
\Ref\Coolidge J.~Coolidge: {\it A treatise on the circle and the sphere\/};
	Oxford University Press, Oxford 1916
\Ref\Guichard C.~Guichard: {\it Sur les syst\`{e}mes triplement
	ind\'{e}termin\'{e}s et sur les syst\`{e}mes triple orthogonaux\/};
	Scientia 25, Gauthier-Villars, Paris 1905
\Ref\Diss U.~Hertrich-Jeromin: {\it \"Uber konform flache Hyperfl\"achen
	in vierdimensionalen Raumformen\/}; PhD thesis, TU Berlin (1994)
\Ref\HertrichJeromin U.~Hertrich-Jeromin: {\it On conformally flat
	hypersurfaces and Guichard's nets\/}; Beitr.\ Alg.\ Geom.\ {\bf 35}
	(1994) 315-331
\Ref\Klein F.~Klein: {\it Vergleichende Betrachtungen \"uber neuere
	geometrische Forschungen (Erlanger Programm)\/}; Math.\ Ann.\
	{\bf 43} (1893) 63-100
\Ref\Palais R.~Palais, C.-L.~Terng: {Critical point theory and submanifold
	geometry\/}; LNM 1353, Springer, New York 1988
\Ref\Salkowski E.~Salkowski: {\it Dreifach orthogonale Fl\"achensysteme\/};
	in Encyclopaedie der mathematischen Wissenschaften III.D 9,
	Teubner, Leipzig 1902
\EndRefs

\Abstract{In the first part, we give a self contained introduction to
the theory of cyclic systems in $n$-dimensional space which can be
considered as immersions into certain Grassmannians. We show how the
(metric) geometries on spaces of constant curvature arise as
subgeometries of M\"obius geometry which provides a slightly new
viewpoint. In the second part we characterize Guichard nets
which are given by cyclic systems as being M\"obius equivalent to
1-parameter families of linear Weingarten surfaces. This provides a
new method to study families of parallel Weingarten surfaces in space
forms. In particular, analogs of Bonnet's theorem on parallel constant
mean curvature surfaces can be easily obtained in this setting.}

\goodbreak\vskip21pt\leftline{\bf Introduction}
	\nobreak\vskip12pt\noindent\kern0pt
Guichard nets were first mentioned by Guichard \Guichard\ who
considered them as a 3-dimensional analog of isothermic nets
\Salkowski. Recently, the first author of the present paper
discovered a close relation to the theory of conformally flat
hypersurfaces in 4-space: any conformally flat hypersurface (in $\R^4$)
carries curvature line coordinates which satisfy the Guichard condition
\HertrichJeromin. In this sense, conformally flat hypersurfaces might
be considered as ``isothermic hypersurfaces\Footnote{It is still not
clear whether the existence of Guichard curvature line coordinates
{\it characterizes} conformally flat hypersurfaces.}''. It seems
remarkable that Guichard already introduces the ``spectral parameter''
which is used in the integrable system approach to conformally flat
hypersurfaces \Diss. In \HertrichJeromin\ it was also shown how
Guichard nets in 3-space --- now considered as special triply
orthogonal systems of surfaces --- do correspond to 3-dimensional
conformally flat hypersurfaces: for certain, very special kinds of
Guichard nets it was possible to characterize the corresponding
conformally flat hypersurfaces geometrically. All of these
``special'' Guichard nets consist of a 1-parameter family of spheres
and two 1-parameter families of channel surfaces, any two surfaces of
different families intersecting orthogonally along a curvature line.
On channel surfaces, one family of curvature lines consists of circles
--- thus, this type of Guichard nets can be considered as cyclic
systems, i.e.\ as given by a 2-parameter family of circles which have
a 1-parameter family of orthogonal surfaces. This was the observation
initiating the present paper: we addressed the problem of classifying
all Guichard nets which come from cyclic systems.

In the first part of this paper we will give a comprehensive
introduction to the theory of cyclic systems. Even though cyclic
systems and M\"obius geometry are very well introduced in the
classical literature (see \Coolidge\ or \Blaschke) it seemed worth to
present this introduction here for two reasons: first, we are going
to use Cartan's method of moving frames which not only allows a
modern formulation of the presented theory but also provides the
tools to characterize the spaces of $m$-dimensional spheres in
$n$-space as certain Grassmannians. Thus, the structure of the spaces
of $m$-spheres becomes very lucid.
And, second, we are going to present a different approach to M\"obius
geometry by introducing it as a supergeometry of the ``metric''
geometries of certain spaces of constant curvature. This new
viewpoint in M\"obius geometry allows to consider geometric problems
in all spaces of constant curvature simultaneously --- as we will
learn in the second part of the present paper, by discussing
Weingarten surfaces and generalizations of Bonnet's theorem on
parallel constant mean curvature surfaces in space forms.
In the concluding section of the first part we will leave the
$n$-dimensional setting and shortly discuss those basic facts in the
theory of triply orthogonal systems in 3-space which we will need for
our discussions on Guichard's nets: systems of three 1-parameter
families of surfaces such that any two surfaces from different
families intersect orthogonally \Salkowski.

In the second part we will present a characterization for Guichard's
nets which come from cyclic systems --- for ``cyclic Guichard nets''.
These turn out to be M\"obius equivalent to 1-parameter families of
parallel Weingarten surfaces in space forms, the Gau\ss{} and mean
curvatures of all surfaces satisfying an affine relation. This is
where the relation of M\"obius geometry and its metric subgeometries
becomes significant. We use this relation to present a way how
various analogs of Bonnet's theorem on parallel constant mean
curvature surfaces in Euclidean space can be obtained (compare
\Palais). Considering a family of parallel linear Weingarten surfaces
in a space form as a cyclic Guichard net the family is naturally
parametrized by an elliptic function. We present relations between
the properties of this elliptic function and the geometry of the
family of parallel Weingarten surfaces --- in particular, we relate
the function's branch points to surfaces of constant mean or
constant Gau\ss{} curvature occuring in the family.

We would like to thank Konrad Voss (ETH Z\"urich) and Hermann Karcher
(University Bonn) for their interest in our work, for fruitful questions
and discussions.

\Part{I: Cyclic Systems}
The main goal of this first part of the present paper is to provide a
modern approach to the theory of cyclic systems --- even though many
of the results presented may be found in the classical literature
(see for example Coolidge's excellent book \Coolidge) it seemed worth
not only to generalize the theory to $n$-dimensional space but also
to point out the relations to symmetric spaces: circles and, more
generally, $m$-spheres can be considered as elements in certain
symmetric spaces. The theory of cyclic systems belongs into the
context of M\"obius geometry. Consequently, we give a short
introduction to M\"obius geometry --- which is slightly different
from the classical approach: instead of considering Euclidean space
or the (metric or conformal) $n$-sphere as the underlying space we
introduce M\"obius geometry on spaces of constant curvature. This
approach allows us to consider all metric geometries simultaneously
as subgeometries of M\"obius geometry --- which will turn out useful
in the second part of the paper. We conclude the first part with
some considerations on triply orthogonal systems --- again, slightly
different from the classical approach (see for example \Salkowski)
since we will work in spaces of constant curvature instead of Euclidean
space.

\Section{The metric geometries}
Later, we are going to introduce M\"obius geometry as a supergeometry
of all the geometries given by the isometry groups on the spaces of
constant curvature $k\in\R$.
For this purpose we define quadric models for the spaces of constant
curvature that will allow us to consider all constant curvature spaces
simultaneously: let $\R^{n+2}_1$ denote the $(n+2)$-dimensional Minkowski
space equipped with a Lorentz scalar product $\langle.,.\rangle$ of signature
$(+,\dots,+,-)$ and $$L^{n+1}:=\{v\in\R^{n+2}_1\,|\,\langle v,v\rangle=0\}$$
its light cone.

\Lemma{~(spaces of constant curvature)}{For any $n_k\in\R^{n+2}_1\setminus
\{0\}$, the intersection $$Q^n_k:=\{v\in L^{n+1}\,|\,\langle v,n_k\rangle=1\}$$
of the light cone with the affine hyperplane $\langle v,n_k\rangle=1$
is a Riemannian space of constant sectional curvature $k=-\langle
n_k,n_k\rangle$.}

In case $k\neq0$ it is immediately clear that $Q_k^n$ has constant
sectional curvature: restricting our attention to the hyperplane
$\langle v,n_k\rangle=1$ we see that $Q_k^n$ is a round sphere of
radius $1\over\sqrt{k}$ or a two-sheeted hyperboloid of ``radius''
$1\over\sqrt{-k}$ in $(n+1)$-dimensional Minkowski space.
In case $k=0$ we choose a point $p_0\in Q_0^n$ and orthogonally
decompose the Minkowski space $\R^{n+2}_1=\R^n\oplus\hbox{span}\{p_0,n_0\}$
--- note that the plane $\hbox{span}\{p_0,n_0\}$ has signature $(+,-)$.
Then, the map $x\in\R^n\mapsto x+p_0-{1\over2}|x|^2n_0\in Q^n_0$ is an
isometry\Footnote{In fact, this isometry is the parametrization often
used in the classical introductions to M\"obius geometry (cf.\ \Blaschke).}.
Thus, for $k\geq0$ we got models for the space forms --- complete and
simply connected --- of curvature $k$ whilst in case $k<0$ we got
spaces which are composed of two copies of a space form of sectional
curvature $k$.

Before discussing the isometry groups of our spaces of constant
curvature we have to understand the relations between these spaces:
consider two quadrics $Q_k^n$ and $Q_{\tilde k}^n$ of constant curvatures
$k$ and $\tilde k$. Then, we can map one onto the other by means of a
rescaling $p\in Q_k^n\mapsto\tilde p=\pm e^up\in Q_{\tilde k}^n$ with a
suitable function $u$ --- note that this map lacks to be defined
for the asymptotic directions $p\perp n_{\tilde k}$, the ``infinity
boundary'', of $Q_{\tilde k}^n$. Since $\langle p,p\rangle\equiv0$ we have
$$\langle d(e^up),d(e^up)\rangle=e^{2u}\langle dp,dp\rangle,$$
i.e.\ these maps between our spaces of constant curvature are conformal.
Moreover, in case of the map $Q_1\to Q_0$ we obtain the classical
stereographic projection after ``correctly'' identifying $Q_1^n\cong
S^n$ and $Q_0^n\cong\R^n$ which convinces us of the following

\Definition{~(generalized stereographic projections)}{The central
projection from one quadric $Q_k^n$ onto another $Q_{\tilde k}^n$ is
called a generalized stereographic projection.}

Clearly, the Lorentz transformations $F\in O_1(n+2)$ which fix $n_k$
act as isometries on $Q_k^n$. Let us try to understand the converse:
given an isometry $f:Q_k^n\to Q_k^n$ we may define a Lorentz transformation
$F:\R^{n+2}_1\to\R^{n+2}_1$ to linearly approximate the isometry
around one point $p\in Q_k^n$: $$Fn_k=n_k,\ Fp=f(p),\ F|_{T_pQ_k^n}=d_pf$$
--- note, that $T_pQ_k^n=\{p,n_k\}^{\perp}$. On space forms, isometries
are uniquely determined by their behaviour at one point. Consequently,
for $k\geq0$, the group of isometries is identical with the group of
those Lorentz transformations that fix $n_k$.
For $k<0$, the situation is a bit more complicated --- in this case
the quadrics $Q_k^n$ are not connected. To learn the specialty of those
transformations coming from Lorentz transformations of the ambient
Minkowski space we project such a hyperbolic quadric $Q_k^n$, $k<0$,
stereographically into a sphere $Q_1$ --- this way we make the
infinity boundary between the two connected components of $Q_k^n$
``visible''. An isometry which comes from a Lorentz transformation
extends smoothly through the infinity boundary\Footnote{On $Q_1$, an
isometry of $Q_k^n$ appears as a conformal transformation --- composed
of an isometry and two stereographic projections.} $\langle
p,n_k\rangle=0$:

\Definition{~(proper isometries)}{The group of proper isometries of
a quadric $Q_k^n$ of constant curvature $k$ is identical with the
group of Lorentz transformations that fix the global normal vector
$n_k$ of $Q_k^n$: $$\hbox{Isom}(Q_k^n)=\{F\in O_1(n+2)\,|\,Fn_k=n_k\}.$$}

For the rest of this paper we will refer to ``proper isometries''
simply as ``isometries''.
Note, that in case $k<0$ the space of proper isometries is still not
the space of hyperbolic motions but a twofold covering of it: there
are proper isometries which interchange the two connected components
of a hyperbolic quadric $Q_k^n$.

\Section{Spheres and Circles}
In order to introduce M\"obius geometry as a supergeometry of all the
metric geometries just discussed we also have to understand the space
of spheres in those spaces $Q_k^n$. To start with, let us consider
hyperspheres as totally umbilic, ``conformally connected\Footnote{Again,
there arise problems with the hyperbolic spaces: here, spheres can be
``connected through the infinity boundary'', i.e.\ a sphere may
consist of two pieces which can be smoothly glued together by adding
a submanifold of the infinity boundary.}'' hypersurfaces in $Q_k^n$.
Given a hypersphere $S$ we may write the umbilicity condition
$dn+hdp=0$ where $p$ varies on $S$, $n$ is a unit normal field to
$S\subset Q_k^n$ and $h$ denotes the (constant) mean curvature of $S$
with respect to $n$. This equation integrates to $n+hp=s$ where
$s\in\R^{n+2}_1$ is a constant vector. Since $n|_p\in T_pQ_k^n=\{p,n_k
\}^{\perp}$ we find that $s\perp p$ is a unit vector and, moreover,
we obtain the sphere's mean curvature as $$h=\langle s,n_k\rangle.$$
At this point, it becomes clear that the vector $s$ in fact {\it
characterizes} the (oriented) sphere $S$: if $p\perp s$ then $s-hp\in
T_pQ_k^n$ and consequently, $p$ is a point on a totally umbilic
hypersurface with (constant) mean curvature $h=\langle s,n_k\rangle$
and unit normal vector $n|_p=s-hp$ in $p$.
Totally umbilic hypersurfaces with vanishing mean curvature are
usually called hyperplanes --- we obtain the hyperplanes in $Q_k^n$
as special hyperspheres: as long as we are interested in M\"obius
geometry rather than in any of its metric subgeometries it is
convenient not to distinguish them.

Since the equation $\langle p,s\rangle=0$ encoding incidence of
a point $p\in Q_k^n$ and a sphere $s\in S^{n+1}_1$ is obviously
independent of the scaling of $p$, spheres are mapped to spheres
by the generalized stereographic projections --- actually, the notion
of a ``sphere'' does not depend on the metric of a space but on its
conformal class only.
As a fundamental invariant in M\"obius geometry also the intersection
angle of two (oriented) spheres can be nicely described in this model:
given two spheres intersecting in a point $p$, $s=n+hp,\tilde s=\tilde
n+\tilde hp\in S^{n+1}_1$, their intersection angle is given by
$$\langle n,\tilde n\rangle=\langle s,\tilde s\rangle.$$

\Lemma{~(oriented spheres)}{The space of oriented hyperspheres in
any of the spaces $Q_k^n$ of constant curvature can be canonically
identified with the Lorentz sphere $S^{n+1}_1\subset\R^{n+2}_1$;
hyperplanes in a $Q_k^n$ may be distiguished by the vanishing of
the mean curvature: $$\langle n_k,s\rangle=0.$$
A point $p\in Q_k^n$ lies on a sphere $s\in S^{n+1}_1$ if and only if
$$\langle p,s\rangle=0.$$
Two spheres $s,\tilde s\in S^{n+1}_1$ intersect orthogonally if and
only if $$\langle s,\tilde s\rangle=0.$$}

From the previous discussions it also becomes clear that two
hyperspheres $s=n|_p+hp$ and $\tilde s=n|_p+\tilde hp=s+(\tilde
h-h)p$ have the same tangent planes in their common point $p$:

\Definition{~(parabolic sphere pencil)}{The 1-parameter family
$s+\R p$ of all oriented hyperspheres sharing tangent planes (and
orientation) in one point $p$ is called a parabolic sphere pencil ---
or an (oriented) hypersurface element.}

Later on, we will not only be interested in hyperspheres but in lower
dimensional spheres, too. Especially, since we are going to discuss
cyclic systems, we will be interested in circles: 1-dimensional spheres.
Similar to the hypersphere case, we may consider $m$-dimensional
spheres in any of the spaces $Q_k^n$ ($m<n$) as totally umbilic
submanifolds of dimension $m$ --- which can be obtained intersecting
$n-m$ totally umbilic (orthogonal) {\it hyper}surfaces: hyperspheres.
Thus, the points $p$ of an $m$-sphere are given by $$p\in L^{n+1}\cap
\{s_1,\dots,s_{n-m}\}^{\perp}$$ where $s_1,\dots,s_{n-m}\in S^{n+1}_1$.
Obviously, this characterization does not depend on the choice of
hyperspheres but only on the $(n-m)$-plane spanned by the $s_i$:

\Theorem{~($m$-spheres)}{The space of (non oriented) $m$-spheres in
any of the spaces $Q_k^n$ of constant curvature can be canonically
identified with the Grassmannian
$$G_+(n-m,m+2)={{O_1(n+2)}\over{O(n-m)\times O_1(m+2)}}$$
of spacelike $(n-m)$-planes of the Minkowski space $\R^{n+2}_1$;
$m$-planes in a $Q_k^n$ are given by those spacelike $(n-m)$-planes
in $\R^{n+2}_1$ which are perpendicular to $n_k$.}

Note that this theorem complements very well our previous statement
on {\it oriented} $(n-1)$-spheres: the Lorentz sphere $S^{n+1}_1$ is
a double cover of the Grassmannian $G_+(1,n+1)$ of spacelike lines in
$\R^{n+2}_1$ --- the two possible orientations of a spacelike line in
$\R^{n+2}_1$ can be interpreted as the orientations of the corresponding
hypersphere. Similarly, $m$-spheres could be oriented by considering the
two possible orientations on spacelike $(n-m)$-planes.

To conclude this section we will discuss 1-spheres --- commonly
called ``circles'' --- a bit more comprehensively: let
$c=\hbox{span}\{s_1,\dots,s_{n-1}\}\in G_+(n-1,3)$ denote a circle.
The points on $c$ are given by lightlike lines\Footnote{Note, that
for the moment we want our arguments to be independent of the ambient
constant curvature spaces $Q_k^n$ --- we want to obtain {\it M\"obius
geometric} notions.} in $c^{\perp}$.
So, we may obtain a parametrization of $c$ by choosing a basis for
$c^{\perp}$: to fit the geometric situation best we choose a {\it
pseudo orthonormal} basis $(s,p,\hat{p})$, i.e.\ $s\in S^{n+1}_1$ is
an (oriented) sphere and $p,\hat{p}\in L^{n+1}$ are two points on $s$
($\langle s,p\rangle=\langle s,\hat{p}\rangle=0$) with $\langle
p,\hat{p}\rangle=1$. Since all the spheres $s_i$ intersect $s$
orthogonally and contain $p$ as well as $\hat{p}$ ($s_i\perp
s,p,\hat{p}$) our circle $c$ is exactly that circle intersecting
the sphere $s$ orthogonally in the two points $p$ and $\hat{p}$:
this yields a description for circles which complements our first
description in some sense.

Given a circle $c$ in this complementary description --- by two
points $p,\hat p$ on a sphere $s$ --- its most general ``arc length''
parametrization into the light cone\Footnote{Again, note that we do
{\it not} require $p_t$ to take values in any of our constant
curvature spaces --- this will turn out convenient in the second
part since we are going to study special coordinate systems there.
To make the parametrization take values in a space of constant
curvature $p_t$ would have to be suitably rescaled.} is given by
$$t\mapsto p_t={1\over g'}(g\,s+p-{1\over2}g^2\hat p)\in L^{n+1}\cap
c^{\perp}\Eqno\Label\ArcLength$$ where $g=g(t)$ denotes {\it any}
function of $t$.
To understand the geometry of this parametrization we fix a projective
scale --- three points --- on the circle: in our setting a somehow
canonical choice would be $p=p_t|_{g=0}$, $\hat{p}=p_t|_{g=\infty}$
and $p_t|_{g=1}=s+p-{1\over2}\hat{p}$. With this choice, $g(t)$ is the
cross ratio\Footnote{In this paper we consider the cross ratio an invariant
of a quadrilateral rather than of a point pair: our cross ratio differs
from the classical one by a transposition of the points.} of the
three scale points and $p_t$: $$g(t)=R(p,p_t,\hat p,p_t|_{g=1})
	=\sqrt{{\langle p,gs+p-{1\over2}g^2\hat p\rangle
	\langle\hat p,s+p-{1\over2}\hat p\rangle}\over
	{\langle gs+p-{1\over2}g^2\hat p,\hat p\rangle
	\langle s+p-{1\over2}\hat p,p\rangle}}.\Eqno\Label\CrossRatio$$
Note, that this cross ratio does not depend on the scaling of the
points in the light cone: it is a conformal invariant, just like
the intersetion angle of two hyperspheres.

An example of a circle which will become important in the second part
is a straight line in any of the constant curvature spaces $Q_k^n$.
Here, we may choose the sphere $s$ in our parametrization as a plane,
i.e.\ $\langle s,n_k\rangle=0$. Since we are considering a straight
line we also have $\langle s_i,n_k\rangle=0$ and, consequently,
$n_k\in\hbox{span}\{p,\hat p\}$. Now, assuming $p$ to be a point in
$Q_k^n$, we find $n_k=-{k\over2}p+\hat p$. For the function $g$, this
yields the differential equation $1=\langle p_t,n_k\rangle={1\over
g'}(1+{k\over4}g^2)$; its solutions --- we fix the initial value
$g(0)=0$, i.e.\ $p_0=p$ ---
$$g(t)=\left\{\matrix{
	{{2\sinh(\sqrt{-k}t)}\over{\sqrt{-k}(1+\cosh(\sqrt{-k}t))}}
		&\hbox{for}&k<0\cr
	t&\hbox{for}&k=0\cr
	{{2\sin(\sqrt{k}t)}\over{\sqrt{k}(1+\cos(\sqrt{k}t))}}
		&\hbox{for}&k>0\cr}\right.$$
provide the arc length parametrizations for straight lines in $Q_k^n$:
$$p_t=\left\{\matrix{
	-{1\over k}\,n_k+{\sinh(\sqrt{-k}t)\,s+\cosh(\sqrt{-k}t)
		({\sqrt{-k}\over2}p-{1\over\sqrt{-k}}\hat{p})\over\sqrt{-k}}
		&\hbox{for}&k<0\cr
	t\,s+p-{1\over2}t^2\hat p\hfil&\hbox{for}&k=0\cr
	-{1\over k}\,n_k+{\sin(\sqrt{k}t)\,s+\cos(\sqrt{k}t)
		({\sqrt{k}\over2}p+{1\over\sqrt{k}}\hat{p})\over\sqrt{k}}
		&\hbox{for}&k>0\cr
	}\right..\Eqno\Label\Geodesic$$
Note, that $p_t$ never reaches the infinity boundary of $Q_k^n$ since
$\langle p_t,n_k\rangle\equiv1\neq0$. In case $k<0$, this means that
$p_t$ only parametrizes half of a straight line, that component which
belongs to $p=p_0$.

At this point we are prepared to introduce

\Section{M\"obius geometry as a supergeometry}
of the metric geometries on the spaces $Q_k^n$ of constant curvature.
To this extend we have to define the group of M\"obius transformations
on $Q_k^n$:

\Definition{~(M\"obius group)}{A transformation\Footnote{Here, we use
the term ``transformation of $Q_k^n$'' in a slightly generalized sense:
we allow a M\"obius transformation to miss the preimage and the image of
the infinity boundary of $Q_k^n$ --- two points for $k=0$ and two spheres
for $k<0$.} of $Q_k^n$ which maps hyperspheres to hyperspheres is called
a M\"obius transformation of $Q_k^n$. The group formed by all M\"obius
tranformations is called the M\"obius group.}

As we noticed earlier, the stereographic projections map spheres in any
$Q_k^n$ to spheres in any $Q_{\tilde k}^n$. Since a Lorentz transformation
$F\in O_1(n+2)$ induces an isometry $F:Q_k^n\to F(Q_k^n)$ it clearly maps
spheres in $Q_k^n$ to spheres in $F(Q_k^n)$, a quadric of the same constant
curvature $k$ but, generally, different from $Q_k^n$. Thus, by composing $F$
with a suitable (unique) stereographic projection we obtain a M\"obius
transformation of $Q_k^n$.
In fact, {\it all} M\"obius transformations of a quadric $Q_k^n$ of
constant curvature can be obtained this way:

\Lemma{~(M\"obius transformations)}{Any M\"obius transformation $\mu$
of $Q_k^n$ is the composition $\mu=\sigma_F\circ F$ of a Lorentz
transformation $F\in O_1(n+2)$ with the (unique) stereographic
projection $\sigma_F:F(Q_k^n)\to Q_k^n$.}

Before attacking the proof of this lemma let us state some facts:
first, it becomes clear that M\"obius transformations are conformal,
i.e.\ they preserve intersection angles between spheres.
And second, the interplay of a M\"obius transformation with the
infinity boundary of $Q_k^n$ in case $k\leq0$ beomes clear: the
stereographic projection $\sigma_F$ lacks to be defined for the
points $p\in Q_k^n$ mapped to the infinity boundary $n_k^{\perp}$
of $Q_k^n$. Similarly, the infinity boundary of $Q_k^n$ will
generally be mapped to a finite region --- to a point for $k=0$
and to a hypersphere in case $k<0$.

In the above construction of a M\"obius transformation, by composing
a Lorentz transformation with a stereographic projection, the (uniquely
determined) stereographic projection was only needed to adjust the
scaling of points correctly. If, for a moment, we identify all the
quadrics of constant curvature by identifying points with lightlike
lines\Footnote{This is the only point where we really adapt the classical
viewpoint: classically, an (absolute) quadric in projective $(n+1)$-space
$\R P^{n+1}$, the ``conformal $n$-sphere'', is considered as the underlying
space for M\"obius geometry --- the Minkowski space $\R^{n+2}_1$ becomes
the space of homogeneous coordinates of $\R P^{n+1}$ the Lorentz product
being fixed (up to scaling) by the conformal $n$-sphere as absolute quadric.
In this model, $m$-spheres are identified with $(n-m-1)$-planes by polarity.},
$p\leftrightarrow\R\cdot v$ where $v\in L^{n+2}$, we can identify a
M\"obius transformation $\mu=\sigma_F\circ F$ with the corresponding
Lorentz transformation $F\in O_1(n+2)$.

Now, that we do not have to care about the ``proper'' scaling of points
any more the proof of our lemma becomes easy: preserving hyperspheres a
M\"obius transformation will also preserve $m$-spheres since those can
be obtained as intersections of hyperspheres. Consequently, a M\"obius
transformation naturally extends to a linear transformation of the Minkowski
space $\R^{n+2}_1$: it maps the spaces $G_+(n-m,m+2)$ of $m$-spheres onto
themselves. But, it is well known that linear transformations of $\R^{n+2}_1$
which preserve the light cone are real multiples of Lorentz transformations.
Thus, given a M\"obius transformation $\mu:Q_k^n\to Q_k^n$, there is a Lorentz
transformation (unique up to sign) $F\in O_1(n+2)$ such that
$\mu=\sigma_F\circ F$ --- this proves the lemma. Moreover,

\Lemma{~(M\"obius group)}{The group of Lorentz transformations is a
(trivial) double covering of the M\"obius group:
$$\hbox{M\"ob}(Q_k^n)\times\{\pm1\}\cong O_1(n+2).$$}

From the way how we introduced the metric geometries it now is
immediately clear that M\"obius geometry is a supergeometry of the
metric geometries\Footnote{In a similar way as M\"obius geometry is
a subgeometry of projective geometry (cf.\ Klein's Erlanger program
\Klein).}:

\Theorem{~(metric subgeometries)}{The geometries of the groups of
motions on the quadrics $Q_k^n$ are subgeometries of the M\"obius
geometry on $Q_k^n$.}

\Section{Envelopes}
Given an immersion $f:M^{n-1}\to Q_k^n$ with unit normal field
$n:M^{n-1}\to S^{n+1}_1$ the normal field may be reinterpreted
as a sphere congruence --- a 2-parameter family of spheres ---
according to our identification of spheres in $Q_k^n$ with unit
vectors in $\R^{n+2}_1$. Since $\langle f,n\rangle=0$ any point
$f(p)$ lies on the corresponding sphere $n(p)$. Moreover, $f$ and
$n(p)$ have first order contact --- $f$ ``touches'' $n(p)$ --- in
$f(p)$ since $\langle d_pf,n(p)\rangle=0$: $n(p)$ can be considered a
common normal of $f(M)$ and $n(p)$ in $f(p)$. Note the different
interpretations for $n(p)$ --- we will use this ambiguity of the
geometric meaning of $n(p)$ repeatedly.
Any other sphere $s\in S^{n+1}_1$ touching $f(M)$ in $f(p)$
lies in one of the parabolic sphere pencils $\pm(n(p)+\R f(p))$,
these two sphere pencils are {\it characterized} by the equations
$\langle s,f(p)\rangle=0$ and $\langle s,d_pf\rangle=0$:

\Definition{~(envelope)}{An immersion $f:M^{n-1}\to Q_k^n$ is said to
envelope a sphere congruence $s:M^{n-1}\to S^{n+1}_1$ if each sphere
$s(p)$ touches $f(M)$ in $f(p)$: $$\matrix{\langle s(p),f(p)\rangle=0
&\hbox{ and }&\langle s(p),d_pf\rangle=0}.\Eqno\Label\Envelope$$}

If we consider the immersion $f$ enveloping a sphere congruence
$s$ to take values in the light cone, $f:M\to L^{n+1}$, rather than
in one of the quadrics $Q_k^n\subset L^{n+1}$ then, the sphere
congruence $s$ can still be interpreted as a (unit) normal field of
$f$ according to \Envelope: $s(p)\in T_{f(p)}L^{n+1}=f(p)^{\perp}$.
On the other hand, since $\langle f,s\rangle\equiv0$ and consequently
$\langle f(p),d_ps\rangle=-\langle s(p),d_pf\rangle$, the immersion
$f$ can as well be interpreted as an (isotropic) normal field of the
sphere congruence $s:M\to S^{n+1}_1$: starting from a hypersphere
congruence $s$ which induces a positive definite metric (i.e.\
$d_ps(T_pM)\in G_+(n-1,3)$) we will be able to find {\it two}
envelopes $f$ and $\hat f$ since the normal bundle of $s$ has
signature $(+,-)$. Geometrically, the points $f(p)$ and $\hat f(p)$
of the two envelopes are given as the intersection points of $s(p)$
and the circles $d_ps(T_pM)$.
---
Focussing on the geometry of a sphere congruence rather than on the
geometry of one of its envelopes it will often be more useful not to
scale the envelopes to lie in a quadric $Q_k^n$ of constant curvature
--- as it might be as well if we are interested in conformal aspects
of the hypersurface's geometry.

To reflect the multiple aspects of a sphere congruence and one of its
envelopes by a more neutral term\Footnote{Note, that in our definition
of a strip we also get rid of any scaling requirements for the point
map $f$ --- as it seems useful when interpreting it as a normal field
of the sphere congruence.} we give the following

\Definition{~(strip)}{A pair $(s,f):M^{n-1}\to S^{n+1}_1\times L^{n+1}$
of smooth maps is called a strip if at least one, $s$ or $f$, is an
immersion with spacelike tangent planes and if
$$\matrix{\langle s,f\rangle=0&\hbox{ and }&\langle s,df\rangle=0}.$$}

To discuss the geometry of strips we will use Cartan's method of moving

\Section{Frames}
Given a strip $(s,f):M\to S^{n+1}_1\times L^{n+1}$ we extend it to a
pseudo orthonormal frame $F:=(s_1,\dots,s_{n-1},s,f,\hat f):M\to
SO_1(n+2)$ into the Lorentz group. The corresponding connection form
$\Phi=F^{-1}dF:TM\to{\fk o}_1(n+2)$ will be of the form
$$\Phi=\left(\matrix{\omega&\eta\cr-\eta^{\ast}&\nu}\right)
	\Eqno\Label\ConnectionForm$$
where $\eta:TM\to\hbox{M}(3\times(n-1))$, $\omega:TM\to{\fk o}(n-1)$ and
$\nu:TM\to{\fk o}_1(3)$ describe the derivative of the circle
congruence $c:=\hbox{span}\{s_1,\dots,s_{n-1}\}$ and the covariant
derivatives on the vector bundles $c$ and $c^{\perp}$, respectively.
This splitting of the connection form $\Phi$ corresponds to the
Cartan decomposition\Footnote{As we will see, circle congruences (and
more generally, congruences of $m$-spheres) can be considered as
immersions into the Grassmannian $G_+(n-1,3)$ ($G_+(n-m,m+2)$).
For that reason, we refer to the Cartan decomposition at this point.}
$${\fk o}_1(n+2)=\left({\fk o}(n-1)\oplus{\fk o}_1(3)\right)\oplus
	\hbox{M}(3\times(n-1))=:{\fk k}\oplus{\fk p}\Eqno\Label\Cartan$$
of the Lie algebra ${\fk o}_1(n+2)$ associated with the symmetric
space $G_+(n-1,3)$ of circles in $Q_k^n$.

First, let us have a closer look at $$\nu=\left(\matrix{0&0&-\nu_s\cr
\nu_s&\nu_f&0\cr0&0&-\nu_f\cr}\right)\Eqno\Label\NU$$ where
$\nu_s=\langle ds,\hat f\rangle$ and $\nu_f=\langle df,\hat f\rangle$.
Since at least one, $s$ or $f$, is assumed to be an immersion we may
choose $s_1,\dots,s_{n-1}$ to span the tangent space of either $s$ or
$f$ --- depending on which one is an immersion and on whether we are
interested in the geometry of $s$ or $f$:

\Definition{~(adapted frames)}{A frame $(s_1,\dots,s_{n-1},s,f,\hat
f)$ of a strip $(s,f)$ is called $s$-adapted (or $f$-adapted) if
$s_1,\dots,s_{n-1}$ span the tangent planes of $s$ (or $f$), i.e.\ if,
in \NU, $\nu_s=0$ (or, $\nu_f=0$).}

In case of an $s$-adapted frame, $\hat f$ describes the second
envelope of the sphere congruence $s$ and the tangent planes
$d_ps(T_pM)$ of $s$ define the congruence of circles orthogonal to
$s$ in its two envelopes $f$ and $\hat f$.
A special case occurs when $s$ is a hyperplane congruence in $Q_k^n$:
then, $\hat f$ is the antipode hypersurface of $f$ ($k>0$), is the
point at infinity ($k=0$) or it is the reflection of $f$ at the
infinity hypersphere $n_k$ ($k<0$).
In case of an $f$-adapted frame, $\hat f$ will usually {\it not} be
the second envelope of $s$ --- however, the circles $d_pf(T_pM)$ will
still intersect the spheres $s(p)$ orthogonally in $f(p)$ and $\hat
f(p)$. If, for example, $f:M\to Q_k^n$ then the circles $d_pf(T_pM)$
are straight lines in $Q_k^n$ since $\langle df,n_k\rangle=0$ and
$\hat f$ will be an envelope of $s$ only if $\hat f$ is a ``parallel
surface'' of $f$, i.e.\ if $\langle s,n_k\rangle=const$.

According to the Cartan decomposition \Cartan\ of ${\fk o}_1(n+2)$ the
Maurer-Cartan equation\Footnote{Here, $[\Phi\wedge\Psi](x,y):=
[\Phi(x),\Psi(y)]-[\Phi(y),\Psi(x)]$ where $\Phi$ and $\Psi$ are
Lie algebra valued 1-forms. In case of a matrix Lie algebra,
$\Phi=(\varphi_{ij})$, (as in our case), where the Lie bracket
becomes the commutator, we may write $$[\Phi\wedge\Phi]=
2\Phi\wedge\Phi:=2\sum_j\varphi_{ij}\wedge\varphi_{jk}.$$}
$d\Phi+{1\over2}[\Phi\wedge\Phi]=0$ splits into the Gau\ss{}-Ricci
equations
$$d\omega+\omega\wedge\omega=\eta\wedge\eta^{\ast}\Eqno\Label\Gauss$$
and
$$d\nu+\nu\wedge\nu=\eta^{\ast}\wedge\eta\Eqno\Label\Ricci$$
and the Codazzi equation
$$d\eta+\omega\wedge\eta+\eta\wedge\nu=0.\Eqno\Label\Codazzi$$
Since we will always work with adapted frames where $\omega$
describes the covariant derivative of $s(M)$ (resp $f(M)$) we will
refer to \Gauss\ as the Gau\ss{} equation and to \Ricci\ as the Ricci
equation. --- From now on, let us assume the frame $F$ to be either
$s$-adapted or $f$-adapted. With the ansatz $$\eta=\left(\matrix{
-\eta_1&\varphi_1&\hat\varphi_1\cr\vdots&\vdots&\vdots\cr
-\eta_{n-1}&\varphi_{n-1}&\hat\varphi_{n-1}\cr}\right)\Eqno\Label\ETA$$
one of the scalar Ricci equations reads $\sum_j\eta_j\wedge\varphi_j=0$
showing that the second fundamental form $\sum_j\eta_j\varphi_j$ of
$s(M)$ with respect to the isotropic normal field $f$ or of $f(M)$
with respect to the unit normal field $s$, respectively, is symmetric.
Consequently, the tangential framing $(s_1,\dots,s_{n-1})$ can be
chosen to simultaneously diagonalize the first and second fundamental
forms.

\Definition{~(curvature framing)}{An ($s$- or $f$-) adapted framing
is called a principal curvature framing of the strip $(s,f)$ if its
tangential part $(s_1,\dots,s_{n-1})$ diagonalizes the second fundamental
form with respect to the induced metric: $$\eta_i\wedge\varphi_i=0.$$
The directions orthogonal to the planes $\eta_i=0$ (or $\varphi_i=0$,
respectively) are then called principal curvature directions of the
strip.}

Passing from an $s$-adapted to the ``nearest'' $f$-adapted principal
curvature framing of the strip $(s,f)$ --- provided both, $s$ as well
as $f$, are immersions --- via $s_i\mapsto s_i+u_if$ with suitable
functions $u_i$, the forms $\eta_i$ and $\varphi_i$ do not change.
Consequently, the principal curvature directions do not change either.

Moreover, a ``conformal deformation'' $(s,f)\mapsto(s+hf,e^uf)$ of a
strip $(s,f)$ yields $\eta_i\mapsto(\eta_i-h\varphi_i)$ and $\varphi_i\mapsto
e^u\varphi_i$ for any corresponding principal curvature framing.
Hence, the principal curvature directions are also not effected by
such a conformal deformation --- passing from a strip $(s,f)$ to the
corresponding immersion $f:M\to Q_k^n$ into a quadric of constant
curvature $k$ with its unit normal field $n:M\to TQ_k^n$ we see that
the principal curvature directions defined above are indeed the principal
curvature directions in the classical sense. We summarize these results in a

\Lemma{}{The principal curvature directions of an immersion $f:M\to
Q_k^n$ are conformally invariant, i.e.\ they coincide with the
principal curvature directions of any strip $(s,e^uf):M\to
S^{n+1}_1\times L^{n+1}$. In particular, the principal curvature
directions are invariant under the stereographic projections
$Q_k^n\to Q_{\tilde k}^n$.}

Given an {\it immersed} sphere congruence $s:M\to S^{n+1}_1$ we have
seen that the two isotropic normal fields in an $s$-adapted framing
$F:M\to O_1(n+2)$ can be interpreted as the two envelopes of $s$.
Assuming both, $f$ as well as $\hat f$, to be immersed the principal
curvature directions of $(s,f)$ and $(s,\hat f)$ will generally {\it
not} coincide. If, however, they do we have $\eta_i\wedge\varphi_i=0$
and $\eta_i\wedge\hat\varphi_i$ at the same time for any $s$-adapted
principal curvature framing $F$. In this case we can arrange to have
$\nu=0$ by possibly rescaling $(f,\hat f)\mapsto(e^uf,e^{-u}\hat f)$
the two isotropic normal fields:
the vector bundle $c^{\perp}=\hbox{span}\{s,f,\hat f\}$ is flat.

\Definition{~(Ribaucour sphere congruence)}{An (immersed) sphere
congruence $s:M\to S^{n+1}_1$ is called a Ribaucour sphere congruence
if its normal bundle is flat: $\eta^{\ast}\wedge\eta=0.$ If its two
envelopes $f$ and $\hat f$ are immersed, this means that their
curvature directions do correspond.}

\Section{Cyclic systems}
Let $F=(s_1,\dots,s_{n-1},s,f,\hat f):M^{n-1}\to O_1(n+2)$ denote an
$s$-adapted framing of a strip $(s,f):M^{n-1}\to S^{n+1}_1\times L^{n+1}$
where $s$ is a Ribaucour sphere congruence and $f,\hat f:M\to L^{n+1}$
are {\it parallel} isotropic normal fields of $s$, i.e.\ $\nu=0$. Then,
$$f_t:={1\over g'(t)}\left(g(t)\cdot s+f-{1\over2}g^2(t)\cdot\hat f\right)
	\Eqno\Label\OrthoSurface$$ will provide simultaneous arc length
parametrizations \ArcLength\ for all circles of the congruence
$c=\hbox{span}\{s_1,\dots,s_{n-1}\}:M\to G_+(n-1,3)$. Moreover, all
circles intersect each hypersurface of the 1-parameter family $(f_t)_t$
orthogonally since $$\langle\D{}{t}f_t,df_t\rangle=0\Eqno\Label\OrthoCondi$$
--- the circle congruence $c$ is what is called a ``cyclic system'':

\Definition{~(cyclic system)}{A circle congruence $c:M^{n-1}\to G_+(n-1,3)$
is called a normal congruence of circles, or a cyclic system, if there
is a 1-parameter family of hypersurfaces (in $Q_k^n$) orthogonal to
all circles.}

Clearly, any immersed sphere congruence defines a circle congruence
which has two orthogonal hypersurfaces --- the two envelopes of the
sphere congruence. Generally, these are the only hypersurfaces which
are orthogonal to all circles of the congruence --- as we will see:
let $c:M^{n-1}\to G_+(n-1,3)$ denote a circle congruence and let
$F=(s_1,\dots,s_{n-1},s,f,\hat f):M\to O_1(n+2)$ be a pseudo
orthonormal framing of $c$. Note, that $F$ is not necessarily a
framing of a strip, i.e.\ $f$ might {\it not} describe an envelope of
the sphere congruence $s$: this fact is responsible for the need to
consider a more general form of the connection form \ConnectionForm\
of $F$: $$\nu=\left(\matrix{0&-\hat\nu_s&-\nu_s\cr\nu_s&\nu_f&0\cr
	\hat\nu_s&0&-\nu_f\cr}\right).$$
With the ansatz \OrthoSurface, $g\equiv1$ --- but $t$ now denoting a
function on $M^{n-1}$, for the orthogonal hypersurfaces of the circle
congruence the orthogonality condition \OrthoCondi\ yields the
following differential equation for $t:M^{n-1}\to\R$:
$$dt={1\over2}t^2\nu_s+t\nu_f+\hat\nu_s.\Eqno\Label\CSDgl$$
The integrability condition for this partial differential equation reads
$$0={1\over2}t^2(d\nu_s+\nu_f\wedge\nu_s)+t(d\nu_f+\hat\nu_s\wedge\nu_s)
	+(d\hat\nu_s+\hat\nu_s\wedge\nu_f)\Eqno\Label\CSInt$$
which, for fixed $p\in M^{n-1}$, is a quadratic polynomial in $t$:

\Theorem{}{If there are more than two hypersurfaces orthogonal to all
circles of a congruence $c:M^{n-1}\to G_+(n-1,3)$, then the circle
congruence is normal.}

Obviously, this statement is due to the fact that a quadratic
polynomial must vanish identically if it has more than two zeros.
Now, the coefficients in \CSInt\ are exactly the coefficients in the
curvature form $d\nu+\nu\wedge\nu$ of the vector bundle $c^{\perp}$
over $M^{n-1}$, i.e.\

\Theorem{}{A circle congruence $c:M^{n-1}\to G_+(n-1,3)$ is normal if
and only if the vector bundle $c^{\perp}$ over $M^{n-1}$ is flat.}

The differential equation \CSDgl\ for $t$ becomes trivial exactly
when the basis fields $s$, $f$ and $\hat f$ of $c^{\perp}$ are
parallel --- then, $t$ can be considered a ``real'' parameter for
the 1-parameter family $(f_t)_t$ of orthogonal hypersurfaces of the
cyclic system. Recalling the geometric interpretation \CrossRatio\
of $g(t)$ in the parametrization \OrthoSurface\ we obtain

\Theorem{}{Any four orthogonal hypersurfaces of a cyclic system
intersect all circles at a fixed cross ratio.}

And finally, since $s:M^{n-1}\to S^{n+1}_1$ is a Ribaucour sphere
congruence, we come back to our starting point:

\Theorem{}{Any two orthogonal hypersurfaces of a cyclic system
envelope a Ribaucour sphere congruence. The circles that intersect
the spheres of a Ribaucour congruence orthogonally in its two envelopes
form a cyclic system.}

\Section{Triply orthogonal systems}
We just learned that any two hypersurfaces orthogonal to the circles
of a cyclic system envelope a Ribaucour sphere congruence.
Consequently, the curvature directions on all orthogonal
hypersurfaces of a cyclic system do correspond: integrating the
$(n-1)$ curvature directions on one orthogonal hypersurface we obtain
$(n-1)$ 1-parameter families of 2-dimensional surfaces --- each
surface built up from the circles along one curvature line --- which
intersect all orthogonal hypersurfaces in their curvature lines.
In case of 3-dimensional ambient space\Footnote{Note, that the
situation is rather special in 3-dimensional ambient space: in higher
dimensions the curvature lines of a hypersurface generally do not
come from a coordinate system.} this yields what is called a ``triply
orthogonal system'' of surfaces:

\Definition{~(triply orthogonal system)}{A system of three 1-parameter
families of surfaces in a 3-dimensional space $Q_k^3$ is called a
triply orthogonal system if any two surfaces from different families
intersect orthogonally.}

Classically, triply orthogonal systems were considered in Euclidean
ambient space (see for example \Salkowski) but the notion of a triply
orthogonal system is obviously conformally invariant: a generalized
stereographic projections $Q_k^3\to Q_{\tilde k}^3$ will map any
triply orthogonal system in $Q_k^3$ onto one in $Q_{\tilde k}^3$
--- and so will any M\"obius transformation do. Since the curvature
directions of a surface in $Q_k^3$ are invariant under the stereographic
projections, too, Dupin's theorem \Salkowski\ on triply orthogonal
systems in Euclidean 3-space holds in spaces of constant curvature:

\Theorem{~(Dupin)}{The surfaces of a triply orthogonal system in
$Q_k^3$ intersect along their curvature lines.}

Applying this theorem to a triply orthogonal system coming from a
cyclic system in $Q_k^3$, we conclude that the two families of
surfaces orthogonal to the orthogonal surfaces of the cyclic system
consist of channel surfaces: these surfaces carry one family of
circular curvature lines and, consequently\Footnote{This is a
consequence of Joachimsthal's theorem.}, each surface envelopes a
1-parameter family of spheres (compare \Coolidge). In fact, this is
the characterization of cyclic systems from the viewpoint of triply
orthogonal systems:

\Theorem{}{A triply orthogonal system comes from a cyclic system if
and only if two of the 1-parameter families of surfaces consist of
channel surfaces.}

Before discussing Guichard's nets in the second part of this paper it
remains to learn some facts about triply orthogonal systems in general:
given a triply orthogonal system in parametric form $(t_1,t_2,t_3)\mapsto
f(t_1,t_2,t_3)\in Q_k^3$, i.e.\ the surfaces of the system being given
by $t_i=const$, we may choose a pseudo orthonormal framing
$F=(n_1,n_2,n_3,f,\hat f):M^3\to O_1(5)$ wherein $n_i$ denote the unit
normal fields --- or, according to our previous identification of the
Lorentz sphere with the space of (oriented) hyperspheres in $Q_k^3$,
the tangent planes --- of the surfaces $t_i=const$ and $\hat f=
n_k+{k\over2}f$ describes the second intersection point of the $n_i$.
The connection form $\Phi=F^{-1}dF$ of such a framing is of the form
$$\Phi=\left(\matrix{
	0&-k_{21}\omega_1+k_{12}\omega_2&-k_{31}\omega_1+k_{13}\omega_3
		&\omega_1&{k\over2}\omega_1\cr
	k_{21}\omega_1-k_{12}\omega_2&0&-k_{32}\omega_2+k_{23}\omega_3
		&\omega_2&{k\over2}\omega_2\cr
	k_{31}\omega_1-k_{13}\omega_3&k_{32}\omega_2-k_{23}\omega_3&0
		&\omega_3&{k\over2}\omega_3\cr
	-{k\over2}\omega_1&-{k\over2}\omega_2&-{k\over2}\omega_3&0&0\cr
	-\omega_1&-\omega_2&-\omega_3&0&0\cr}\right)\Eqno\Label\TOSCF$$
where $\omega_i=l_idt_i$ with Lam\'{e}'s functions $l_i:=|\D{f}{t_i}|$
and $k_{ij}=-{1\over l_il_j}\D{}{t_i}l_j$ give the principal curvatures
of the surfaces $t_i=const$ in $t_j$-direction. The Maurer-Cartan
equation $d\Phi+\Phi\wedge\Phi=0$ reduces to Lam\'{e}'s equations
$$\matrix{
\matrix{k=e_1k_{12}+e_2k_{21}-k_{12}^2-k_{21}^2-k_{31}k_{32}\cr
	k=e_2k_{23}+e_3k_{32}-k_{23}^2-k_{32}^2-k_{12}k_{13}\cr
	k=e_3k_{31}+e_1k_{13}-k_{31}^2-k_{13}^2-k_{23}k_{21}\cr}\cr
		\vbox to 2ex{}\cr
\matrix{0=e_1k_{23}+k_{13}(k_{21}-k_{23})\cr
	0=e_2k_{31}+k_{21}(k_{32}-k_{31})\cr
	0=e_3k_{12}+k_{32}(k_{13}-k_{12})\cr}\hskip1.5em
\matrix{0=e_1k_{32}+k_{12}(k_{31}-k_{32})\cr
	0=e_2k_{13}+k_{23}(k_{12}-k_{13})\cr
	0=e_3k_{21}+k_{31}(k_{23}-k_{21})\cr}\cr
}\Eqno\Label\Lame$$
where $e_i:={1\over l_i}\D{}{t_i}$ are the unit vector fields in
$t_i$-direction. In fact, as the Maurer-Cartan equation for the
connection form \TOSCF\ Lam\'{e}'s equations are exactly the
conditions on three functions $l_i$, $i=1,2,3$, to determine a triply
orthogonal system:

\Theorem{~(Lam\'{e})}{Three functions $l_1$, $l_2$ and $l_3$ are the
Lam\'{e} functions of a triply orthogonal system in a quadric $Q_k^3$
of constant curvature $k$ if and only if they satisfy Lam\'{e}'s
equations \Lame\ with $e_i={1\over l_i}\D{}{t_i}$ and
$k_i=-{1\over l_il_j}\D{}{t_i}l_j$.}

In some situations --- especially when examining Guichard's nets ---
it is more convenient to allow the parametrization $f$ {\it not} to
take values in one quadric $Q_k^3$ but, more generally, in the light
cone $L^{n+1}$. Then, the unit vectors $n_i/\!/\D{f}{t_i}$ will not
longer describe the tangent plane congruences of the surfaces
$t_i=const$ but sphere congruences which are enveloped by the
surfaces. On the connection form \TOSCF\ this has the effect that
$d\hat f=\sum_{i,j=1}^3b_{ij}n_i\omega_j$ where the functions
$b_{ij}$ can be determined from the Maurer-Cartan equation:
$$\matrix{
\matrix{b_{11}+b_{22}=e_1k_{12}+e_2k_{21}-k_{12}^2-k_{21}^2-k_{31}k_{32}\cr
	b_{22}+b_{33}=e_2k_{23}+e_3k_{32}-k_{23}^2-k_{32}^2-k_{12}k_{13}\cr
	b_{33}+b_{11}=e_3k_{31}+e_1k_{13}-k_{31}^2-k_{13}^2-k_{23}k_{21}\cr}\cr
		\vbox to 2ex{}\cr
\matrix{b_{21}=e_1k_{23}+k_{13}(k_{21}-k_{23})\cr
	b_{32}=e_2k_{31}+k_{21}(k_{32}-k_{31})\cr
	b_{13}=e_3k_{12}+k_{32}(k_{13}-k_{12})\cr}\hskip1.5em
\matrix{b_{31}=e_1k_{32}+k_{12}(k_{31}-k_{32})\cr
	b_{12}=e_2k_{13}+k_{23}(k_{12}-k_{13})\cr
	b_{23}=e_3k_{21}+k_{31}(k_{23}-k_{21})\cr}\cr
}\Eqno\Label\BIJ$$
with $e_i$ and $k_{ij}$ as above. Note, that $b_{ij}=b_{ji}$. The
integrability conditions for Lam\'{e}'s functions now become third
order differential equations --- the conformal flatness\Footnote{In
fact, Lam\'{e}'s equations \Lame\ are the conditions on the metric
$\sum_{i=1}^3l_i^2dt_i^2$ to have constant curvature $k$.} of the
metric $\langle df,df\rangle=\sum_{i=1}^3l_i^2dt_i^2$:
$$\matrix{
\matrix{e_1b_{23}-k_{23}b_{31}-k_{13}b_{23}&=&
	e_3b_{21}-k_{21}b_{13}-k_{31}b_{21}\cr
	e_2b_{31}-k_{31}b_{12}-k_{21}b_{31}&=&
	e_1b_{32}-k_{32}b_{21}-k_{12}b_{32}\cr
	e_3b_{12}-k_{13}b_{23}-k_{32}b_{12}&=&
	e_2b_{13}-k_{13}b_{32}-k_{23}b_{13}\cr}\cr
		\vbox to 2ex{}\cr
\matrix{0=e_1b_{22}-e_2b_{21}+k_{12}(b_{11}-b_{22})
	+k_{21}(b_{12}+b_{21})+k_{32}b_{31}\cr
	0=e_2b_{33}-e_3b_{32}+k_{23}(b_{22}-b_{33})
	+k_{32}(b_{23}+b_{32})+k_{13}b_{12}\cr
	0=e_3b_{11}-e_1b_{13}+k_{31}(b_{33}-b_{11})
	+k_{13}(b_{31}+b_{13})+k_{21}b_{23}\cr}\cr
		\vbox to 2ex{}\cr
\matrix{0=e_2b_{11}-e_1b_{12}+k_{21}(b_{22}-b_{11})
	+k_{12}(b_{21}+b_{12})+k_{31}b_{32}\cr
	0=e_3b_{22}-e_2b_{23}+k_{32}(b_{33}-b_{22})
	+k_{23}(b_{32}+b_{23})+k_{12}b_{13}\cr
	0=e_1b_{33}-e_3b_{31}+k_{13}(b_{11}-b_{33})
	+k_{31}(b_{13}+b_{31})+k_{23}b_{21}\cr}\cr
}\Eqno\Label\GenLame$$
where $e_i={1\over l_i}\D{}{t_i}$, $k_{ij}=-{1\over l_il_j}\D{}{t_i}l_j$
and $b_{ij}$ are defined by \BIJ. These are the equations we will use
later --- in place of the original Lam\'{e} equations \Lame\ --- when
discussing Guichard's nets in $Q_k^3$.

\Part{II: Guichard's nets}
In the second part we are going to classify those Guichard nets ---
special triply orthogonal systems --- which come from cyclic systems.
A Guichard net can be considered as a 3-dimensional analog of an
isothermic net in the plane \Guichard --- the surfaces of a Guichard net
divide the ambient space into infinitesimal rectangular parallelepipeds
(any two surfaces of different families intersect orthogonally along
curvature lines) such that two of the six diagonal rectangles are squares:

\Definition{~(Guichard net)}{A triply orthogonal system is called a
Guichard net if, with a suitable choice of $\varepsilon_i\in\{1,i\}$,
its Lam\'{e} functions $l_i$ satisfy $$\sum_{i=1}^3(\varepsilon_il_i)^2
=0.\Eqno\Label\GuichardCondi$$}

As we will see, a cyclic system which gives rise to a Guichard net
consists of a family of parallel Weingarten surfaces in some space of
constant curvature\Footnote{This behaviour is similar to that of
isothermic Willmore surfaces (Thomsen's theorem \Blaschke): a surface
belonging to two M\"obius geometric surface classes turns out to
belong to a metric surface class --- minimal surfaces in spaces of
constant curvature.}. Thus, as a byproduct, we obtain several
generalizations of Bonnet's theorem on parallel surfaces of constant
mean curvature in Euclidean space \Palais.
To complete the discussion, we will try to give an ``estimate'' for
the generality of those Guichard nets coming from cyclic systems by
discussing the effect of the assumption to come from a cyclic system
on Lam\'{e}'s equations \GenLame.

\Section{Cyclic systems}
Let us start with a cyclic system $c:M^2\to G_+(2,3)$ --- we assume
that the family of its orthogonal surfaces is given by $f_t:M^2\to L^4$
such that $t$ is the arc length on all circles simoultanously, i.e.\
$|\D{}{t}f_t|\equiv1$ which can always be achieved by a suitable
scaling of $f_t$ into the light cone. Denoting $f=f_0$ we find that
$$f_t={1\over g'(t)}[g(t)\,s+f-{1\over2}g^2(t)\,\hat f]\Eqno\Label\OF$$
with a Ribaucour sphere congruence $s$, its second envelope $\hat f$
and where $g$ denotes some function with $g(0)=0$ and $g'(0)=1$ (see
\OrthoSurface).
The condition \GuichardCondi\ for the corresponding triply orthogonal
system to be a Guichard net reads --- since all circles are simultaneously
parametrized by arc length ---
$$\left|\D{}{t_1}f_t\right|^2+\varepsilon^2\left|\D{}{t_2}f_t\right|^2=1$$
where $(t_1,t_2)$ are suitable curvature line coordinates and
$\varepsilon\in\{1,i\}$, depending on the position of the diagonal
infinitesimal squares in the Guichard net relativ to the circle
direction. The Guichard condition for $t=0$ gives us
$$df=\left\{\matrix{\cos u\cdot s_1dt_1+\sin u\cdot s_2dt_2\cr
	\cosh u\cdot s_1dt_1+\sinh u\cdot s_2dt_2\cr}\right.$$
with $s_1,s_2:M^2\to S_1^4$ and a suitable function $u:M^2\to\R$.
With the ansatz
$$\matrix{
	ds&=&\left\{\matrix{
	-(a_1\cos u-a_2\sin u)s_1dt_1-(a_1\sin u+a_2\cos u)s_2dt_2\cr
	-(a_1\cosh u+a_2\sinh u)s_1dt_1-(a_1\sinh u+a_2\cosh u)s_2dt_2\cr
	}\right.\hfill\cr
	d\hat f&=&\left\{\matrix{
	(b_1\cos u-b_2\sin u)s_1dt_1+(b_1\sin u+b_2\cos u)s_2dt_2\cr
	(b_1\cosh u+b_2\sinh u)s_1dt_1+(b_1\sinh u+b_2\cosh u)s_2dt_2\cr
	}\right.\hfill\cr}\Eqno\Label\Ansatz$$
where $a_i,b_i:M^2\to\R$ denote suitable functions the Guichard
condition becomes
$$\matrix{g'^2&=&
	(1-a_1g-{1\over2}b_1g^2)^2+\varepsilon^2(a_2g+{1\over2}b_2g^2)^2
	\hfill\cr&=&
	[1-(a_1+\varepsilon a_2)g-{1\over2}(b_1+\varepsilon b_2)g^2]\cdot
	[1-(a_1-\varepsilon a_2)g-{1\over2}(b_1-\varepsilon b_2)g^2].\cr}$$
This equation has some interesting consequences: first,
the function $g$ has to be an elliptic function. Since it does not
depend on $(t_1,t_2)$ its branch points do not either. Hence, $a_i$
and $b_i$ are constant and, consequently, there is a constant
vector $n_k:=b_2s+(a_1b_2-a_2b_1)f+a_2\hat f$ perpendicular to all
$s_1$ and $s_2$: $$dn_k=d[b_2s+(a_1b_2-a_2b_1)f+a_2\hat f]=0.$$
Without loss of generality\Footnote{The vector $n_k$ vanishes if and
only if $a_2=b_2=0$. But in this case all surfaces $f_t$ are totally
umbilic --- they form a sphere pencil. The corresponding conformally
flat hypersurfaces were completely classified in \HertrichJeromin.},
we may assume $n_k\neq0$: then, the cyclic system $c=\hbox{span}
\{s_1,s_2\}:M^2\to G_+(2,3)$ consists of straight lines in the quadric
$Q_k^3$ of constant sectional curvature $k=-\langle n_k,n_k\rangle$
corresponding to the vector $n_k$, and the surfaces $f_t$ are
parallel in $Q_k^3$.

\Proposition{Any cyclic system which defines a Guichard net is a
normal line congruence in some quadric $Q_k^3$ of constant curvature.}

Finally, we find that the Maurer-Cartan equation for the adapted
principal curvature framing $F=(s_1,s_2,s,f,\hat f):M^2\to O_1(5)$
reduces to some version of the sine-Gordon equation:
$$\matrix{0&=&(\D{}{t_1}\D{}{t_1}-\D{}{t_2}\D{}{t_2})u
	+{1\over2}([a_1^2-a_2^2+2b_1]\sin2u + [a_1a_2+b_2]\cos2u)\,,\hfill\cr
	0&=&(\D{}{t_1}\D{}{t_1}+\D{}{t_2}\D{}{t_2})u
	+{1\over2}([a_1^2+a_2^2+2b_1]\sinh2u +[a_1a_2+b_2]\cosh2u)\cr}$$
depending on whether $\varepsilon=1$ or $\varepsilon=i$, respectively.

\Section{Parallel Weingarten surfaces}
Now, that we know that our cyclic system is in fact a normal line
congruence in some space of constant curvature we may choose a more
adapted parametrization for the family of orthogonal surfaces:
after choosing a ``base surface'' $f=f_0$ out of the family of
parallel surfaces --- and assuming $f$ to actually take values in
$Q_k^3$, i.e.\ $\langle f,n_k\rangle=1$ --- we might fix the Ribaucour
sphere congruence $s$ to be its tangent plane congruence, i.e.\
$\langle s,n_k\rangle=0$. For the vector $n_k$ defining the quadric
$Q_k^3$ this means $$n_k=-{k\over2}f+\hat f$$ and the second envelope
$\hat f$ of $s$ becomes trivial: it is the antipode surface of $f$,
the point at infinity or the reflection of $f$ at the infinity
boundary depending on whether $k>0$, $k=0$ or $k<0$, respectively. In
\Ansatz, this is reflected by the fact that $2b_1=k$ and $b_2=0$.
This way, we have fixed the frame $(s,f,\hat f)$ which we use to
parametrize the circles of the congruence.

Instead of using arc length parametrizations \Geodesic\ --- which
would lead to unpleasant calculations --- we choose the easiest
parametrization possible for the tangent plane congruences $s_t$ of
the parallel surfaces $f_t$:
$$s_t:={1\over\sqrt{1+kt^2}}\left(s+t\cdot[{k\over2}f+\hat f]\right).$$
Asking now all $f_t:M^2\to Q_k^3$ to take values in $Q_k^3$, i.e.\
fixing the scaling of $f_t$ ``correctly'', gives us
$$\matrix{f_t&=&{1\over\sqrt{1+kt^2}}\left(f-t\cdot[s
		+{t\over1+\sqrt{1+kt^2}}\cdot n_k]\right),\hfill\cr
	\hat f_t&=&{1\over\sqrt{1+kt^2}}\left(\hat f-{k\,t\over2}\cdot[s
		+{t\over1+\sqrt{1+kt^2}}\cdot n_k]\right).\cr
	}\Eqno\Label\Parametrization$$
Herein, the range where $t$ is running is restricted by the condition
$1+kt^2>0$: in case $k<0$ this prevents us from running through the
infinity boundary --- as $t\to\pm{1\over\sqrt{-k}}$ the points of
$f_t$ approach the intersection points of each circle with the
infinity sphere $n_k$, as the points of $\hat f_t$ do ``from the
other side''. In case of Euclidean ambient space, $k=0$, $f_t$
parametrizes each circle up to the point $n_0\equiv\hat f_t$ at
infinity which is approached as $t\to\pm\infty$. Finally, in case
$k>0$, $f_t$ also parametrizes just half of each great circle,
$\hat f_t$ taking the other half; the two antipode points
$f_{\pm\infty}=\mp{1\over\sqrt{k}}s+{1\over k}n_k$ which lie
symmetric with respect to $f_0$ and $\hat f_0$ are never reached
but just approached by $f_t$ and $\hat f_t$ as $t\to\pm\infty$.

Since $s_t$ is the tangent plane congruence of $f_t$, i.e.\ it can be
interpreted as the normal field of $f_t:M^2\to Q_k^3$, we can easily
calculate the first and second fundamental forms of $f_t$ --- there
is no need to calulate those of $\hat f_t$ since $\hat f_t$ is {\it
isometric} to $f_t$ in $Q_k^3$: $$\matrix{I_t=\langle df_t,df_t\rangle&
	\hbox{and}&I\!I_t=-\langle ds_t,df_t\rangle}.$$
With \Ansatz\ --- remember that we have $b_1={k\over2}$ and $b_2=0$
since $s$ is the tangent plane congruence of $f=f_0$ in $Q_k^3$ ---
we find $$\matrix{k_1&=&\left\{\matrix{
	{{(a_1-kt)\cos u-a_2\sin u}\over{(1+a_1t)\cos u-a_2t\sin u}}\hfill\cr
	{{(a_1-kt)\cosh u+a_2\sinh u}\over{(1+a_1t)\cosh u+a_2t\sinh u}}\cr
		}\right.&\hbox{and}&k_2&=&\left\{\matrix{
	{{(a_1-kt)\sin u+a_2\cos u}\over{(1+a_1t)\sin u+a_2t\cos u}}\hfill\cr
	{{(a_1-kt)\sinh u+a_2\cosh u}\over{(1+a_1t)\sinh u+a_2t\cosh u}}\cr
		}\right.}$$
for the principal curvatures of $f_t$. Since in both cases, $\varepsilon=1$
and $\varepsilon=i$, both principal curvatures are given in terms of
{\it one} function, $u$, the surfaces $f_t$ clearly are Weingarten.
Moreover, it is easy to see that the (extrinsic) Gau\ss{} curvature
$k_t=k_1k_2$ and the mean curvature $H_t={1\over2}(k_1+k_2)$ of the
surfaces $f_t:M^2\to Q_k^3$ satisfy an affine relation
$0=c_K(t)\cdot K_t+2c_H(t)\cdot H_t+c(t)$ where $$\matrix{
c_K(t)&=&[a_1^2+(\varepsilon a_2)^2]\cdot t^2+2a_1\cdot t+1,\hfill\cr
c_H(t)&=&a_1k\cdot t^2+[k-(a_1^2+(\varepsilon a_2)^2)]\cdot t-a_1,\hfill\cr
c(t)&=&k^2\cdot t^2-2a_1k\cdot t+[a_1^2+(\varepsilon a_2)^2].\hfill\cr}
	\Eqno\Label\Cs$$
Note, that the sign of $$c_K(t)\cdot c(t)-c_H^2(t)=\varepsilon^2a_2^2
(1+kt^2)^2\Eqno\Label\CaseInvariant$$ is an invariant of the family
of parallel Weingarten surfaces: it determines the position of the
infinitesimal squares in the Guichard net relative to the family's
normal direction.

\Proposition{The orthogonal surfaces of a cyclic system defining
a Guichard net are (parallel) linear Weingarten surfaces in a quadric
$Q_k^3$ of constant curvature, i.e.\ their (extrinsic) Gau\ss{} and
mean curvature satisfy an affine relation
$$c_K(t)\cdot K_t+2c_H(t)\cdot H_t+c(t)=0.$$}

Since we lately changed the parametrization \Parametrization\ of
the cyclic system we need to derive the condition on the triply
orthogonal net to form a Guichard net again: to that extend suppose
$t=t(r)$ is a function of a new parameter $r$. Then, the condition
$|\D{f_t}{r}|^2=|\D{f_t}{t_1}|^2+\varepsilon^2|\D{f_t}{t_2}|^2$
to form a Guichard net becomes
$$t'^2=[1+kt^2]\cdot c_K(t)\Eqno\Label\Elliptic$$
which, again, is the equation of an elliptic function. From this
equation we first see that in case $k<0$ of hyperbolic ambient space
Guichard nets do {\it not} extend through the infinity sphere: for
$t\to\pm{1\over\sqrt{-k}}$ the Guichard net becomes singular. Thus,
any Guichard net defined through a cyclic system is a normal line
congruence in a {\it space form} rather than just in one of the
constant curvature quadrics $Q_k^3$. Also, any cyclic Guichard net
becomes singular at the zeros of $c_K$: these are surfaces of
constant mean curvature which might occur in the family of parallel
Weingarten surfaces.

\Proposition{A Guichard net given by a cyclic system becomes singular
at any constant mean curvature surface present in the family of
orthogonal surfaces of the cyclic system as well as it becomes
singular at the infinity sphere in case of hyperbolic ambient space.}

Examining the occurance of ``special'' surfaces --- characterized
through the zeros of $c_K$, $c_H$ or $c$ --- in a family of parallel
Weingarten surfaces is interesting on its own; but, at the moment, we
would like to postpone this topic and rather discuss the

\Section{Construction of Guichard nets}
which come from cyclic systems: in the preceeding two paragraphs we
learned that those Guichard nets always come from a family of
parallel Weingarten surfaces whose Gau\ss{} and mean curvature
satisfy an affine relation. To understand the reverse, we start with
a Weingarten surface $f:M^2\to Q_k^3$ whose principal curvatures
$k_1$ and $k_2$ satisfy $0=c_Kk_1k_2+c_H(k_1+k_2)+c$ and study the
triply orthogonal system given by the family of its parallel surfaces:
since surfaces of constant mean curvature are singular for cyclic
Guichard nets we exclude them --- hence, we may assume $c_K=1$.
Now, the principal curvature coordinates of $f$ can be fixed in a
canonical way: setting $a_1:=-c_H$ and $a_2:=\sqrt{|c-c_H^2|}$ and
$u$ denoting a suitable function, we make the ansatz $$\matrix{
	k_1=\left\{\matrix{a_1-a_2\tan u\hfill\cr a_1+a_2\tanh u\cr}\right.
		&\hbox{and}&
	k_2=\left\{\matrix{a_1+a_2\cot u\hfill\cr a_1+a_2\coth u\cr}\right.
	}\Eqno\Label\PCAnsatz$$
in case $c-c_H^2>0$ and in case $c-c_H^2<0$, respectively --- at this
point we exclude $c-c_H^2=0$: surfaces with $0=(k_1+c_H)(k_2+c_H)$,
i.e.\ which have one constant principal curvature. Besides for sphere
pieces, surfaces of this kind do {\it not} occur as orthogonal surfaces
of a cyclic system which defines a Guichard net (compare \CaseInvariant).
With the above ansatz for the principal curvatures the Codazzi equations
for $f$ read in principal curvature coordinates $(t_1,t_2)$
$$\matrix{\matrix{
	\D{}{t_2}\left({\sqrt{\langle\D{}{t_1}f,\D{}{t_1}f\rangle}
		\over\cos u}\right)=0\hfill\cr
	\D{}{t_2}\left({\sqrt{\langle\D{}{t_1}f,\D{}{t_1}f\rangle}
		\over\cosh u}\right)=0\hfill\cr
			}&\hbox{and}&\matrix{
	\D{}{t_1}\left({\sqrt{\langle\D{}{t_2}f,\D{}{t_2}f\rangle}
		\over\sin u}\right)=0\hfill\cr
	\D{}{t_1}\left({\sqrt{\langle\D{}{t_2}f,\D{}{t_2}f\rangle}
		\over\sinh u}\right)=0\hfill\cr}}$$
showing that with $\varepsilon\in\{1,i\}$ --- $\varepsilon=1$ if
$c-c_H^2>0$ and $\varepsilon=i$ if $c-c_H^2<0$ --- we can assume
$$\left|\D{}{t_1}f\right|^2+\varepsilon^2\left|\D{}{t_2}f\right|^2\equiv
	const.$$
Now, parametrizing the family of parallel surfaces of our Weingarten
surface --- as for example in \Parametrization\ with the unit normal
field $s=n:M^2\to S_1^4$ of $f$ --- we see that running through it
with the ``correct speed'' \Elliptic\ provides us with a Guichard net.

\Theorem{~(cyclic Guichard nets)\PLabel\MT}{The orthogonal surfaces of
any cyclic system which defines a Guichard net are parallel linear
Weingarten surfaces in a space form, i.e.\ their Gau\ss{} an mean
curvatures satisfy an affine relation $$c_KK+2c_HH+c=0$$ where
$c_K\neq0$ as well as $c_Kc-c_H^2\neq0$ (no surfaces of constant
mean curvature or with a constant principal curvature occur%
\Footnote{Besides the degenerate case of sphere pieces: those
Guichard nets were discussed in \HertrichJeromin.}).
Conversely, the normal line congruence of such a Weingarten surface
always defines a cyclic Guichard net.}

To estimate the ``amount of generality'' of those

\Section{Guichard's nets}
coming from cyclic systems we will study the effect of the assumption
to come from a cyclic system Lam\'{e}'s equations \GenLame\ for any
Guichard net: choosing an appropriate scaling for the parametrization
$f:M^3\to L^4$ of a Guichard net in a quadric $Q_k^3$ of constant
curvature we may assume $l_3\equiv1$ for one of its Lam\'{e} functions.
Then, the condition to be Guichard reads $l_1^2+\varepsilon^2l_2^2=1$
which leads to the ansatz $l_1=\cos(\varepsilon w)$ and
$l_2={1\over\varepsilon}\sin(\varepsilon w)$ where $w:M^3\to\R$
denotes a suitable function. Introducing functions
$$\matrix{w_0&:=&(\D{w}{t_1})^2+\varepsilon^2(\D{w}{t_2})^2+(\D{w}{t_3})^2
		\hfill\cr
	w_1&:=&-\varepsilon\D{}{t_1}\D{w}{t_3}\cdot\cot(\varepsilon w)
		\hfill\cr
	w_2&:=&\varepsilon\D{}{t_2}\D{w}{t_3}\cdot\tan(\varepsilon w)
		\hfill\cr
	w_3&:=&{\varepsilon\over\sin(2\varepsilon w)}\left(
		(\D{}{t_1}\D{}{t_1}-\varepsilon^2\D{}{t_2}\D{}{t_2})w
		-\D{}{t_3}\D{w}{t_3}\cdot\cos(2\varepsilon w)
		\right)\hfill\cr}$$
the generalized Lam\'{e}'s equations reduce to four differential equations
$$\matrix{\D{w_i}{t_j}=\D{w_j}{t_i}&\hbox{and}&
	\D{w_1}{t_1}+\varepsilon^2\D{w_2}{t_2}+\D{w_3}{t_3}=
	\varepsilon^2\D{w_0}{t_3},}\Eqno\Label\GuiLame$$
i.e.\ $(w_1,w_2,w_3)$ is a gradient and $\hbox{div}(w_1,w_2,w_3)=
\D{}{t_3}|\hbox{grad}w|^2$. Comparing the present ansatz to our first
ansatz \OF\ for a Guichard net coming from a cyclic system and using
\Ansatz\ we see that in case the Guichard net comes from a cyclic
system we have
$$\matrix{\cos(\varepsilon w)&=&{1\over g'}(
	[1-a_1g-{1\over2}b_1g^2]\cos(\varepsilon u)+
	[a_2g+{1\over2}b_2g^2]\varepsilon\sin(\varepsilon u)
		)\hfill\cr
	\sin(\varepsilon w)&=&{\varepsilon\over g'}(
	[1-a_1g-{1\over2}b_1g^2]{1\over\varepsilon}\sin(\varepsilon u)-
	[a_2g+{1\over2}b_2g^2]\cos(\varepsilon u)
		)\hfill\cr}$$
and, consequently, the function $w$ has to split: $w(t_1,t_2,t_3)=
u(t_1,t_2)+v(t_3).$ On Lam\'{e}'s equations \GuiLame\ this has a
radical effect: it remains only {\it one} single equation
$$\matrix{{\varepsilon\over\sin(2\varepsilon w)}\left(
	(\D{}{t_1}\D{}{t_1}-\varepsilon^2\D{}{t_2}\D{}{t_2})u-
	v''\cos(2\varepsilon w)\right)-\varepsilon^2(v')^2
	=c=const}$$
which is equivalent to
$$\matrix{(\D{}{t_1}\D{}{t_1}-\varepsilon^2\D{}{t_2}\D{}{t_2})u&=&
	[(c+\varepsilon^2v'^2)\cos(2\varepsilon v)
	-v''\varepsilon\sin(2\varepsilon v)]\cdot{1\over\varepsilon}
		\sin(2\varepsilon u)\cr&+&
	[(c+\varepsilon^2v'^2){1\over\varepsilon}\sin(2\varepsilon v)
	-v''\cos(2\varepsilon v)]\cdot
		\cos(2\varepsilon u).\cr}$$
Herein, the coefficients of $\sin(2\varepsilon u)$ and $\cos(2\varepsilon u)$
have to be constant which ``splits'' the equation: for the function
$u=u(t_1,t_2)$ we obtain a version of the sine-Gordon (sinh-Gordon)
equation and for the function $v=v(t_3)$ we get a modified pendulum
equation $$c+\varepsilon^2v'^2=r_0\cos(2\varepsilon[v-v_0]).$$
Thus, as in our previous ansatzes, the circles of a cyclic Guichard
net are simultaneously parametrized by an elliptic function $v$.
These discussions may convince the reader that those Guichard nets
in a quadric $Q_k^3$ coming from cyclic systems are quite special ---
even though we are far from being able to solve the system \GuiLame\
it seems reasonable to expect the existence of more general solutions
$w$ that do not split into a function $u$ of two and one, $v$, of one
variable only\dots

Because of the relation between systems of parallel Weingarten
surfaces in space forms and cyclic Guichard nets --- which provides
us with a M\"obius geometric characterization for those systems of
parallel Weingarten surfaces --- we may use the M\"obius geometric
setting to study linear Weingarten surfaces.
To present this principle we will give simple proofs for various
generalizations of

\Section{Bonnet's theorem}
on the existence of a parallel constant mean curvature surface to a
given surface of constant mean curvature \Palais.
From our main theorem on cyclic Guichard nets \MT\ we know that
the family of parallel surfaces of a linear Weingarten surface ---
which has non constant principal curvatures --- in a space form of
sectional curvature $k$ are linear Weingarten surfaces, i.e.\ the
(extrinsic) Gau\ss{} and mean curvature of all surfaces $f_t$ in the
family satisfy an affine relation $$0=c_K(t)\cdot K_t+2c_H(t)\cdot
H_t+c(t).$$ In this setting, surfaces of constant Gau\ss{} curvature,
of constant mean curvature or with constant sum ${1\over k_1}+{1\over k_2}$
of the curvature radii are characterized by the zeros of $c_H$, $c_K$
and $c$, respectively. If we parametrize the family $f_t$ as in
\Parametrization\ then $c_K$, $c_H$ and $c$ become the quadratic
polynomials \Cs\ with $a_1=-c_H(0)$ and $a_2=\sqrt{|c(0)-c_H^2(0)|}$
--- here, we assume that $f_0$ is {\it not} a surface of constant
mean curvature\Footnote{This is no restriction: to a surface of
constant mean curvature in any space form there exist plenty of
parallel surfaces which do not have constant mean curvature.} so that
the ansatz \PCAnsatz\ for the principal curvatures works and that we
can assume $c_K(0)=1$, without loss of generality.
Thus, generalizations of Bonnet's theorem can be obtained by studying
the zeros of quadratic polynomials (compare \Palais) --- the only
thing which remains unpleasant is the calculation of the distances
between two surfaces of the family: to that purpose we have to
integrate $$\int|\D{}{t}f_t|dt=\int{dt\over1+kt^2}.$$
Let us summarize these observations in a

\Theorem{~(meta theorem)}{In the M\"obius geometric setting, the
study of parallel surfaces of constant Gau\ss{} and mean curvature
and with constant sum of the curvature radii reduces to the study of
zeros of quadratic polynomials.}

In case $\varepsilon=1$, i.e.\ $c_Kc-c_H^2>0$, neither $c_K$ nor $c$
have real zeros showing that there are no constant mean curvature
surfaces or surfaces with constant sum of their curvature radii in
the family. $c_H$, on the other hand, has always zeros --- in case
$k<0$, one of them lies in $(-{1\over\sqrt{-k}},{1\over\sqrt{-k}})$.
Consequently, there is always (at least) one surface of constant
Gau\ss{} curvature present which we can choose as the ``base
surface'' $f_0$ of the family. But $c_H(0)=0$ means that $c_H$ is,
in fact, linear and hence has exactly {\it one} zero --- provided it
does not vanish identically. In case of Euclidean an hyperbolic
ambient spaces, this means that there is exactly {\it one} surface of
constant Gau\ss{} curvature in the family. In case of elliptic
ambient space, it is easy to see that the surfaces $f_{\pm\infty}$
have constant Gau\ss{} curvature $$\matrix{K_{\pm\infty}=\lim_{t\to\infty}
{-1\over c_K(t)}(c(t)+2c_H(t)H_t)=-{k^2\over a_2^2}},$$ too.
Thus, in case of elliptic ambient space there exist {\it four}
surfaces of constant Gau\ss{} curvature in the family which have
distances ${1\over2\sqrt{k}}\pi$, i.e.\ they divide each normal great
circle in quarters. Since the antipode surface of a surface clearly
has the same curvatures as the surface itself there occur two values
$K_0$ and $K_{\infty}$ for the constant Gau\ss{} curvatures --- these
satisfy $K_0K_{\infty}=k^2$.

In case $\varepsilon=i$, i.e.\ $c_Kc-c_H^2<0$, the situation becomes
more interesting:
if $k>0$, i.e.\ in case of elliptic ambient space, all of the
functions $c_K$, $c_H$ and $c$ have real zeros: as for $\varepsilon=1$,
there are two parallel surfaces of constant Gau\ss{} curvatures $K_1$
and $K_2$, $K_1K_2=k^2$, in distance $\pi\over2\sqrt{k}$ --- and their
antipode surfaces --- there are two (antipode) pairs of parallel surfaces
of constant mean curvature $\pm{1\over2}(\sqrt{K_1}-\sqrt{K_2})$ in
distance $d={1\over\sqrt{k}}\hbox{arctan}\sqrt{k\over K_1}$ from the
$K_1$-surfaces and two pairs of constant $\pm{1\over2k}(\sqrt{K_1}
-\sqrt{K_2})$ sum of their curvature radii, in the same distance $d$
from the $K_2$-surfaces.

If $k=0$, i.e.\ in case of Euclidean ambient space, the function $c_H$
becomes linear, $c_K$ has two real zeros and $c$ has no zeros: we
obtain Bonnet's classical theorem --- provided $c_H$ does not vanish
identically in which case we are left with a family of surfaces of
constant sum of their curvature radii, parallel to a minimal surface.

If $k<0$, i.e in case of hyperbolic ambient space, we observe the widest
variety of cases: if $c_H$ does not vanish identically there can either
be one ore no surface of constant Gau\ss{} curvature in the family.
To a surface of constant Gau\ss{} curvature there exist either two
parallel surfaces of constant mean curvature or two parallel surfaces
with constant sum of their curvature radii\Footnote{It seems remarkable
that in case of elliptic ambient space the {\it extrinsic} Gau\ss{}
curvature of a surface has to be a positive constant in order to have
parallel constant mean curvature surfaces while, in hyperbolic ambient
space, this condition meets the {\it intrinsic} Gau\ss{} curvature.}.
If no surface of constant Gau\ss{} curvature is present in the family
there is either one surface of constant mean curvature, one surface with
constant sum of the curvature radii or one surface of either type.

Viewing the family of parallel linear Weingarten surfaces as a
Guichard net corresponds to the ``correct'' choice of the family
parameter --- $t$ becomes an elliptic function \Elliptic\ of a new
parameter $r$. We already learned that the branch points of $t$ have
a geometric meaning (if they are real) for the family of Weingarten
surfaces: they describe the infinity sphere and surfaces of constant
mean curvature. The cross ratio of the branch points of $t$ turns out
to be real or to lie on the unit circle depending on whether
$\varepsilon\sqrt{k}$ is real or imaginary; consequently, the
underlying torus --- on which the elliptic function is defined when
viewed as a function of a complex variable --- is a rectangular or a
rhombic torus. Thus, if the type $\varepsilon$ of the Weingarten
surfaces a Guichard net is built of is known the type of the
underlying torus corresponds to the type of the ambient space ---
the ambiguous case of cross ratio 1 (where the torus degenerates to a
cylinder) corresponding to Euclidean space.
Moreover, we can establish various relations between the cross ratio
of the branch points of the elliptic function $t$ and geometric
quantities arising from the family of parallel Weingarten surfaces:
for example, in case of cross ratio $-1$ (square torus) the zeros of
$c_K$ and $c$ coincide, i.e.\ all surfaces of constant mean curvature
present in the family are minimal, and the (extrinsic) Gau\ss{}
curvature of any surface of constant Gau\ss{} curvature equals the
ambient space's curvature.

We leave the complete analysis of the situation to the interested reader.

\References

\bye